\documentclass[a4paper]{amsart}

\usepackage{array}

\usepackage{amsmath}                        
\usepackage{amssymb}
\usepackage[english]{babel}                 %
\usepackage[english]{varioref}              
\usepackage{amsfonts}                       
\usepackage{fancyhdr}
\usepackage{enumerate}
\usepackage{amsthm}
\usepackage{multirow}

\newcounter{satz}

\newtheorem{cor}[satz]{Corollary}
\newtheorem{lem}[satz]{Lemma}
\newtheorem{prop}[satz]{Proposition}

\theoremstyle{definition}
\newtheorem*{rem}{Remark}
\newtheorem*{ex}{Example}

\newcommand{\LieG}{{\mathrm G}}
\newcommand{\LieK}{{\mathrm K}}

\newcommand{\LieSpin}{{\mathrm{Spin}}}
\newcommand{\LieO}{{\mathrm{O}}}
\newcommand{\LieZ}{{\mathrm{Z}}}
\newcommand{\LieSO}{{\mathrm{SO}}}

\newcommand{\LieSp}{{\mathrm{Sp}}}
\newcommand{\LieSU}{{\mathrm{SU}}}
\newcommand{\LieU}{{\mathrm{U}}}

\newcommand{\lieG}{{\mathfrak g}}
\newcommand{\lieA}{{\mathfrak{a}}}
\newcommand{\lieB}{{\mathfrak b}}
\newcommand{\lieC}{{\mathfrak{c}}}
\newcommand{\lieD}{{\mathfrak d}}
\newcommand{\lieK}{{\mathfrak k}}
\newcommand{\lieP}{{\mathfrak p}}

\newcommand{\lieO}{{\mathfrak o}}

\newcommand{\lieSpin}{{\mathfrak{spin}}}

\newcommand{\Id}{{\mathrm I}}
\newcommand{\Ad}{{\mathrm{Ad}}}
\newcommand{\ad}{{\mathrm{ad}}}

\newcommand{\A}{{\mathcal{A}}}
\newcommand{\Exp}{{\mathrm{Exp}}}

\def\R{{\mathbb R}}
\def\N{{\mathbb N}}
\def\Z{{\mathbb Z}}
\def\C{{\mathbb C}}
\def\H{{\mathbb H}}

\begin{document}
\title{`Spindles' in symmetric spaces}

    \author[Peter Quast]{Peter Quast}
    \address{D\'{e}partement de math\'{e}matiques, Universit\'e de Fribourg, 1700 Fribourg, Switzerland}
    \curraddr{Institut f\"{u}r Mathematik, Universit\"{a}t Augsburg, 86135 Augsburg, Germany}
    \email{peter.quast@unifr.ch,\; peter.quast@math.uni-augsburg.de}
    \thanks{Author supported by the Swiss National Science Foundation Grant PBFR-106367}
    \date{\today}
    \dedicatory{}

\begin{abstract}
    We study families of submanifolds in symmetric spaces of compact
    type arising as exponential images of $s$-orbits of variable radii.
    If the $s$-orbit is symmetric such submanifolds are the most
    important examples of adapted submanifolds, i.e.\ of submanifolds of
    symmetric spaces with curvature invariant tangent and normal spaces.
\end{abstract}

 \maketitle


\section{Introduction}
Among Riemannian manifolds the symmetric spaces form an important
class of examples. On the one hand they are sufficiently general to
serve as examples for a lot of geometric phenomena. On the other
hand they enjoy enough algebraic structure to make explicit
calculations possible. The most important algebraic structure on a
symmetric space is the Lie triple product on its tangent spaces
given by the curvature tensor. Thus it is natural to ask which
submanifolds of a given Riemannian symmetric space have curvature
invariant tangent and normal spaces at each point. Such a
submanifold is called {\it adapted}. Naitoh classified locally these
submanifolds using the theory of Grassmann geometries (cf.\
\cite{Nait-98} and the references given there). It turns out that if
the ambient symmetric space $M$ is irreducible and of rank greater
than 1 then a full adapted submanifold, i.e. one not contained in a
proper totally geodesic subspace (also called a {\it normal curved
Lie triple} (cf.~\cite{Esch-02})), is locally isomorphic to a
symmetric $s$-orbit of $M.$ An {\it $s$-orbit} is a connected
component of an isotropy orbit of a symmetric space. A result of
Ferus (cf.\ \cite{Feru-80}) shows that symmetric $s$-orbits, also
known as symmetric {\it $R$-spaces}, are exactly the closed
extrinsically symmetric submanifolds of euclidean space. A
submanifold $S$ of $M$ is called {\it extrinsically symmetric} if
for each point $q\in S$ there exists an isometry $f_q$ of $M$
satisfying $f_q(q)=q,\; f_q(S)=S$ and $\mathrm{D}f_q(q)X=-X$ if $X$
is tangent to $S$ and $\mathrm{D}f_q(q)X=X$ if $X$ is normal to $S.$
The isometry $f_q$ is called the {\it extrinsic symmetry} of $S$ at
$q,$ e.g.~if $M$ is a euclidean space then $f_q$ is the reflection
along the normal space of $S$ at  $q$

In this work we study more closely families of exponential images of
$s$-orbits, in particular the cases where the $s$-orbits are
symmetric. These families form {\it chains of spindles} in the
ambient space. The number of knots in a chain of spindles is called
the {\it spindle number.} The submanifolds occurring in such chains
of spindles of symmetric $s$-orbits are described in \cite[pp.\ 263,
264]{Be-Co-Ol-03}, but the description given there is not complete.

Given a Lie triple $(\lieP,R)$ of compact type (i.e.\ a vector space
$\lieP$ equipped with a triple product $R$ satisfying the algebraic
curvature identities with the additional property that $R(X,Y)$ is a
derivation of $R$ for all $X,Y\in\lieP$ and such that the orthogonal
symmetric Lie algebra obtained by the Cartan construction is of
compact type (cf.\ \cite{Esch-02}, \cite{Helg-78})), there are two
canonical symmetric spaces associated with it. On the one hand the
associated simply connected space which is the universal Riemannian
covering of all symmetric spaces associated with $(\lieP,R).$ On the
other hand the {\it adjoint space} which is covered by any symmetric
space associated with $(\lieP,R).$ We show that spindle numbers of
adjoint spaces are 1 (see section \ref{section-adjoint}) and
describe the relation between the spindle number and the center of
the isometry group (see section \ref{section-isometry}). Further we
study the extrinsic geometry of  the chains of spindles associated
with symmetric $s$-orbits (see section \ref{section-extrinsic}) and
determine the spindle numbers associated with symmetric $s$-orbits
of the classical simply connected symmetric spaces of compact type
(see section \ref{section-examples}). It turns out that any natural
number can be realized as the spindle number associated with a
Veronese type embedding of some projective space.

We refer to \cite{Helg-78} for the general theory on symmetric
spaces and to \cite{Be-Co-Ol-03} for a detailed description  of
$s$-orbits.

It is my pleasure to thank J.-H.\ Eschenburg and many other
colleagues from Augsburg University for many valuable discussions
and helpful remarks. Moreover I am  most grateful to E.\ A.\ Ruh for
his support during the past years. This work was supported by the
Swiss National Science Foundation Grant PBFR-106367.

\section{Chains of spindles and spindle
numbers}\label{section_Chains-of-spindle}

Let $M$ be a Riemannian symmetric space of compact type and let $p$
be some point in $M.$ Denote by $\LieG$ its isometry group and by
$\LieK$ the isotropy group of $p.$ The geodesic symmetry of $M$ at
$p$ will be denoted by $s_{p}$ and the corresponding involuting
automorphism $\Ad(s_p)$ of the Lie algebra $\lieG $ of $\LieG$ by
$\sigma.$ The Cartan decomposition of $\lieG$ associated with these
data is $\lieG=\lieK\oplus\lieP,$ where the $(+1)$-eigenspace
$\lieK$ of $\sigma$ is the Lie algebra of $\LieK$ and the
$(-1)$-eigenspace $\lieP$ of $\sigma$ is identified with $T_{p}M$
(cf.\ \cite{Helg-78}).

Let $\gamma$ be a closed geodesic on $M$ emanating from $p$ and let
$\xi\in\lieP$ be its tangent vector at $p.$ Considering a maximal
torus (flat) in $M$ containing $\gamma$ one sees that, after some
suitable change of the parametrization of $\gamma,$ the imaginary
parts of the (purely imaginary) eigenvalues of the anti-selfadjoint
homomorphism $\ad(\xi)$ are integers. If these integers are
relatively prime we call the parametrization of $\gamma$ and the
(non-zero) element $\xi\in\lieP$ {\it canonical.} In this article
$\xi$ always denotes a canonical element of $\lieP$ and $i$ always
the imaginary unit $\sqrt{-1}.$ It is well known that $i\nu$ is an
eigenvalue of $\ad(\xi)$ if and only if $-i\nu$ is an eigenvalue of
$\ad(\xi).$

A point $q$ in $M$ is called an {\it antipode} of $p$ if $q$ is the
midpoint of a closed geodesic emanating from $p.$ The set $\A(p)$ of
antipodes of $p$  is the fixed point set of $s_p.$ Isolated points
in $\A(p)$ are called {\it poles} of $p$ and connected components of
$\A(p)$ having positive dimension are called {\it polars} of $p.$
Polars are examples of adapted submanifolds: Let $q$ be a point of
some polar $M^+(q)$ of $p.$ On the one hand $M^+(q)$ is totally
geodesic in $M.$ On the other hand the connected component of the
fixed point set $s_q\circ s_p$ containing $q,$ called a {\it
meridian}, is a totally geodesic submanifold whose tangent space at
$q$ is the normal space of $M^+(q)$ (cf.~\cite{Ch-Na-78}). Any
connected component of $\A(p)$ is an orbit of $\LieK_0,$ the
identity component of $\LieK,$ in $M$ (cf.~\cite{Ch-Na-78}). Hence
let $q$ be a point of some polar $M^+(q)$ of $p$ and let $\gamma$ be
the closed geodesic emanating from $p$ whose midpoint is
$q=\gamma(t_0)$ then $M^+(q)=\Exp_p(t_0\cdot M^{\xi}),$ where
$M^{\xi}=\Ad_{\LieG}(\LieK_0)\xi\subset\lieP$ is the $s$-orbit of
$\xi=\gamma'(0)$ and $\Exp_p$ the Riemannian exponential map of $M$
at $p.$ The objects introduced in this paragraph were extensively
studied by Nagano and Nagano and Tanaka in a series of subsequent
papers with the same title (see the references given in
\cite{Na-Ta-00}).

Let now $\xi$ be a canonical element in $\lieP$ and let $M^{\xi}$ be
its $s$-orbit. The family of submanifolds
$\left(M^{\xi}_t\right)_{t\in\R},\; M^{\xi}_t=\Exp_{p}\left(t\cdot
M^{\xi}\right)$ is called the {\it chain of spindles} in $M$
associated with the canonical element $\xi.$ A submanifold
$M^{\xi}_t$ which is a point is called a {\it knot}, a family of
submanifolds between two consecutive knots is called a {\it spindle}
and a submanifold in the middle between two consecutive knots is
called a {\it centriole} (cf.~\cite{Ch-Na-88}) of the chain of
spindles. Notice that the set of the geodesics symmetries at the
knots acts transitively on the set of spindles of a given chain of
spindles. Let $t_{\xi}$ be the period of the closed geodesic
$\gamma_{\xi}$ defined by $\xi,$ i.e.\
$\gamma_{\xi}(t_{\xi})=\gamma_{\xi}(0)=p$ and $p\notin
\gamma_{\xi}\left((0,t_{\xi})\right).$ The {\it spindle number}
$\lambda(M,\xi)$ associated with the pair $(M,\xi)$ is the number of
knots of the chain of spindles associated with $\xi$ the geodesic
$\gamma_{\xi}(t)$ meets in the interval $[0,t_{\xi}).$ Since
$p=\gamma_{\xi}(0)$ is a knot, the spindle number is at least 1.
Observe that $\lambda(M,\Ad_{\LieG}(k)\xi)=\lambda(M,\xi)$ for
$k\in\LieK.$

\begin{ex}
The best known example of a chain of spindles is formed by the
`parallels of latitude' on a sphere: Let $M$ be the $2$-dimensional
standard sphere $S^2,\; p$ its north pole and $\xi$ any non-zero
vector in $T_pS^2.$ The corresponding $s$-orbit is the one
dimensional sphere $S^1$ in $T_pS^2.$ The spindle defined by $\xi$
is just the family of parallels of latitude. The knots are the north
and the south pole and the spindle number is 2.
\end{ex}

Given a canonical element $\xi$ of $\lieP,$ we consider the
decomposition of $\lieG,$
    $$\lieG=\lieG_{-\nu_r}\oplus\lieG_{-\nu_{r-1}}\oplus\;\hdots\;\oplus\lieG_{-\nu_{1}}\oplus
    \lieG_0\oplus\lieG_{\nu_{1}}\oplus\;\hdots\;\oplus\lieG_{\nu_{r-1}}\oplus\lieG_{\nu_r},$$
where $\lieG_{\nu}$ is the $(i\nu)$-eigenspace of $\ad({\xi}).$ Then
$\lieK$ and $\lieP$ can be decomposed as follows:
$$\begin{array}{llllll}
    \lieK=\lieK_+\oplus\lieK_-, & \text{where} & \lieK_-=\sum\limits_{j=1}^r\lieK_{\nu_j}, &
    \lieK_{\nu_j}=\lieK\cap\left(\lieG_{-\nu_j}\oplus\lieG_{\nu_{j}}\right) & \text{and} &
    \lieK_+=\lieK\cap\lieG_0;\\
    \lieP=\lieP_+\oplus\lieP_-, & \text{where} & \lieP_-=\sum\limits_{j=1}^r\lieP_{\nu_j}, &
    \lieP_{\nu_j}=\lieP\cap\left(\lieG_{-\nu_j}\oplus\lieG_{\nu_{j}}\right) & \text{and} &
    \lieP_+=\lieP\cap\lieG_0.\\
\end{array}$$
Observe that $\lieP_-\cong T_{\xi}M^{\xi}$ and $\lieP_+\cong
N_{\xi}M^{\xi}.$ Let $X'$ be an element of $\lieK_-$ and let
$X=\ad(X')\xi=\sum_{j=1}^rX_{\nu_j}\in\lieP_-,\;
X_{\nu_j}\in\lieP_{\nu_j}.$ Consider the geodesic variation
    $$\alpha:\R\times\R\longrightarrow M,\quad
    (s,t)\longmapsto\gamma_{\Ad_{\LieG}(\exp(sX'))\xi}(t).$$
Its variation vector field for $s=0$ is the Jacobi field $J_X$ along
$\gamma_{\xi}$ given by
\begin{equation}\tag{$\ast$}\label{Jacobi-field}
    J_X(t)=\sum\limits_{j=1}^r{1\over{\nu_j}}\sin(\nu_jt)X_{\nu_j}(t),
\end{equation}
where $X_{\nu_j}(t)$ is the parallel vector field along
$\gamma_{\xi}$ with $X_{\nu_j}(0)=X_{\nu_j}.$ Assume
$X_{\nu_j}(t)\neq 0$ for all $j\in\{1,\hdots,r\}.$ Since the vector
fields $X_{\nu_j}(t)$ are linearly independent and the eigenvalues
$\{\nu_j\}$ relatively prime, $J_X$ vanishes if and only if
$t\in\Z\pi.$ This shows that $M^{\xi}_t$ is a knot if and only if
$t\in\Z\pi$ and a centriole if and only if $t\in\left(\Z+{1\over
2}\right)\pi.$ Moreover the length of $\gamma_{\xi}$ is
$\lambda(M,\xi)\cdot\pi\cdot |\xi|.$

Two submanifolds $S_1$ and $S_2$ of $M$ are called {\it
extrinsically isometric,} if there exists an isometry $f$ of $M$
mapping $S_1$ onto $S_2,$ i.e.\ $f(S_1)=S_2.$

\begin{lem}\label{spindle=symmetry}
    The submanifolds $M^{\xi}_{n\pi+t}$ and
    $M^{\xi}_{n\pi-t},\; n\in\Z$ are extrinsically isometric.
\end{lem}
\begin{proof}
The geodesic symmetry at the knot $M^{\xi}_{n\pi}$ maps
$M^{\xi}_{n\pi+t}$ onto $M^{\xi}_{n\pi-t}.$
\end{proof}

\subsection{Spindle numbers of adjoint
spaces}\label{section-adjoint} Denote by $\overline{\LieG}$ the
adjoint group $\Ad(\LieG)$
 of $\LieG.$ For an element $g\in\LieG$ we
denote by $\overline{g}$ the corresponding element in
$\overline{\LieG}.$ Consider further the subgroup $\overline{\LieK}$
of $\overline{\LieG}$ formed by the elements of $\overline{\LieG}$
commuting with $\overline{s_p}.$ Then
$\overline{M}=\overline{\LieG}/\overline{\LieK}$ is a symmetric
space, called the {\it adjoint space} of the orthogonal symmetric
Lie algebra $(\lieG,\sigma).$ Notice that the isotropy actions of
$\overline{\LieK}_0$ and $\LieK_0$ on $\lieP$ have the same orbits.
The adjoint space $\overline{M}$ is covered by all symmetric spaces
associated with $(\lieG,\sigma)$ (cf.\ \cite[p.\ 327]{Helg-78}). We
denote by $\pi$ the covering map
$M\stackrel{\pi}{\longrightarrow}\overline{M}.$

\begin{prop}\label{spindle_number_adjoint_space}
  Spindle numbers of
  adjoint spaces are 1.
\end{prop}

\begin{proof}
    Since $\xi$ is canonical, $\overline{\exp}(\pi\xi) =\Ad(\exp(\pi\xi))$ is an
    element of order 2, i.e.\ $\overline{\exp}(\pi\xi)^2=\overline{\exp}(2\pi\xi)=\Id.$
    Since $s_p\circ
    \exp(\pi\xi)\circ s_p=\exp(-\pi\xi),$ we observe that
    $\Ad(\exp(\pi\xi)\circ s_p)= \Ad(s_p\circ\exp(\pi\xi))$. Hence
    $\overline{\exp}(\pi\xi)$ is an element of $\overline{\LieK}$ and
    $\lambda(\overline{M},\xi)$ is 1.
\end{proof}

Proposition \ref{spindle_number_adjoint_space} shows that
$\lambda(M,\xi)$ describes how many times the geodesic in
$\overline{M}$ defined by $\xi$ is covered by the corresponding
geodesic in $M.$ More precisely $\pi$ maps the knots of the spindle
chain $\left(M^{\xi}_t\right)_{t\in\R}$ in $M$ onto $\pi(p)$ and the
centrioles onto the polar (which is not a pole, since the spindle
number is 1) of the geodesic in $\overline{M}$ defined by $\xi.$
This shows that centrioles of chains of spindles are totally
geodesic in $M$ (cf.~\cite{Ch-Na-88}).

\subsection{Spindle numbers and centers of isometry
groups}\label{section-isometry}

\begin{prop}\label{Spindle_center}
    The spindle number $\lambda({M},\xi)$ divides
    the double of the order of the center $\LieZ(\LieG_0)$ of the identity component $\LieG_0$ of $\LieG,$
    i.e.
    $$ \lambda(M,\xi)\;\mid\; 2\cdot \mathrm{card}\left(\LieZ(\LieG_0)\right).$$
    If $\lambda(M,\xi)$ is odd then  $\lambda(M,\xi)$ already divides
    the cardinality of $\LieZ(\LieG_0).$
\end{prop}
\begin{proof} Let $g=\exp(\pi\xi)\in\LieG_0.$ Since $\xi$ is canonical, we have
        $\Ad(g^2)=e^{\ad(2\pi\xi)}=\Id.$
        Thus $g^2$ lies in $\LieZ(\LieG_0).$ Since
        $g^{\lambda(M,\xi)}=\exp(\lambda(M,\xi)\pi\xi)$ is an element of $\LieK,$ we see that
        $g^{2\lambda(M,\xi)}$ lies in $\LieZ(\LieG_0)\cap\LieK,$ i.e.\
        $g^{2\lambda(M,\xi)}=\Id.$ Let $\mu$ be the smallest positive integer
        satisfying $g^{2\mu}=\Id.$ Then $\mu$ divides the cardinality of $\LieZ(\LieG_0)$
        and $\lambda(M,\xi)$ is a multiple of $\mu,$
        i.e.\ $\lambda(M,\xi)=n\mu$ for some positive integer $n.$ Since $p=g^{2\mu}p=\exp(2\pi\mu\xi)p,$
        we have $2\mu\geq \lambda(M,\xi)$  and hence $2\geq n.$ If $n=1$ then $\lambda(M,\xi)$ divides already the
        cardinality of $\LieZ(\LieG_0).$ If $n=2$ then $\lambda(M,\xi)$ divides
        $2\cdot \mathrm{card}\left(\LieZ(\LieG_0)\right).$
\end{proof}

The next example shows that Proposition \ref{Spindle_center} is
optimal: Let $M=\LieSU(6)/\LieSO(6)$ and $\xi=i\left(
\begin{array}{cc}
  -{5\over 6} & 0 \\
  0 & {1\over 6}\Id_5 \\
\end{array}%
\right).$  Then the spindle number
$\lambda(\LieSU(6)/\LieSO(6),\xi)$ is $6$ (see section
\ref{section-examples}) and the cardinality of the center of
$\LieG_0\cong{\LieSU(6)/\Z_2}$ is $3.$

\subsection{Spindle numbers of products}
Let $M_1$ and $M_2$ be two symmetric spaces of compact type and let
$M=M_1\times M_2.$ Let $p=(p_1,p_2)$ be a point in $M.$ By $\LieG_i$
we denote the isometry group of $M_i$ and by $\LieK_i$ the isotropy
group of $p_i,\; i=1,2.$ The Cartan decomposition of the Lie algebra
of the isometry group $\LieG$ of $M$ is
$\lieG=\lieG_1\oplus\lieG_2=\lieK\oplus \lieP,$ where
$\lieK=\lieK_1\oplus\lieK_2$ is the Lie algebra of the isotropy
group $\LieK$ of $p$ and $\lieP=\lieP_1\oplus\lieP_2\cong T_pM$
where $\lieP_i\cong T_{p_i}M_i,\; i=1,2.$ Assume $\xi_i\in
\lieP_i,\; i=1,2$ to be canonical. Then $\xi=(\xi_1,\xi_2)\in\lieP$
is canonical and the spindle number $\lambda(M,\xi)$ is the smallest
common multiple of the spindle numbers $\lambda(M_1,\xi_1)$ and
$\lambda(M_2,\xi_2).$


\section{Chains of spindles of extrinsically symmetric
type}\label{section-extrinsic}

Among the canonical elements there might be some whose $s$-orbit
 $M^{\xi}$ is extrinsically symmetric in the euclidean
space $\lieP.$ This happens exactly if
$$
    \ad(\xi)^3=-\ad(\xi),
$$
i.e.\ if the eigenvalues of $\ad(\xi)$ are $\pm i$ and $0$ (cf.\
\cite{Feru-80}). Such elements will be called of {\it extrinsically
symmetric type}. In this section $\xi$ always denotes such an
element. Considering again the geodesic variation field
(\ref{Jacobi-field}) one sees that chains of spindles associated
with elements of extrinsically symmetric type have the following
characterizing property: All submanifolds $M_t^{\xi}$ which are not
knots have the same dimension.

If $\xi$ is of extrinsically symmetric type the Jacobi field given
in equation (\ref{Jacobi-field}) simplifies:
$$J_X(t)={\mathrm{D}\Exp_p}\mid_{t\xi}(tX)=\sin(t)X(t).$$
Since parallel transports along curves are linear isometries, this
shows that the tangent space
$T_{\gamma_{\Ad_{\LieG}(k)\xi}(t)}M^{\xi}_t, \, t\notin\pi\Z$ is
just the parallel transport of $T_{\Ad_{\LieG}(k)\xi}M^{\xi}\cong
\Ad_{\LieG}(k)\lieP_-$ along $\gamma_{\Ad_{\LieG}(k)\xi}.$ Take now
an element $Y\in\lieP_+\cong N_{\xi}(M^{\xi}).$ The Jacobi field
$J_Y(t)={\mathrm{D}\Exp_p}\mid_{t\xi}(tY)$ along $\gamma_{\xi}$ is
$J_Y(t)=tY(t),$ where $Y(t)$ is the parallel vector field along
$\gamma_{\xi}$ defined by $Y.$ Hence the normal space
$N_{\gamma_{\Ad_{\LieG}(k)\xi}(t)}M^{\xi}_t, \, t\notin\pi\Z$ is the
parallel transport of $N_{\Ad_{\LieG}(k)\xi}M^{\xi}\cong
\Ad_{\LieG}(k)\lieP_+$ along $\gamma_{\Ad_{\LieG}(k)\xi},\;
k\in\LieK_0.$ Since $\lieP_+$ and $\lieP_-$ are curvature invariant
if $\xi$ is of extrinsically symmetric type, and since the curvature
tensor of $M$ is parallel, we get: {\it The submanifolds
$M^{\xi}_{t}, \, t\notin\pi\Z$ are adapted.} In fact, these
submanifolds form the most prominent class of adapted submanifolds
in symmetric spaces of compact type (cf.\ \cite{Esch-02},
\cite{Nait-98}).

A submanifold $S$ of $M$ is called {\it reflective} if there exists
an isometry $r$ of $M$ fixing $S$ pointwise and reversing the
geodesics of $M$ which intersect $S$ perpendicularly or,
equivalently, if $S$ is a connected component of the fixed point set
of an involuting isometry of $M,$ e.g.\ polars and meridians are
reflective. The involuting isometry $r$ is called the {\it
reflection} of $M$ in $S.$ A reflective submanifold in a symmetric
space is totally geodesic and moreover extrinsically symmetric: Let
$r$ be the reflection of $M$ in $S$ and let $q$ be any point in $S.$
Then $r\circ s_q$ is the extrinsic symmetry of $S$ at $q.$
Conversely, a totally geodesic extrinsically symmetric submanifold
is reflective.

\begin{prop}{\bf (cf.~\cite[pp.\ 256, 257]{Be-Co-Ol-03},
\cite{Nait-86})}\label{extr-sym-adapted-reflective}
    Let $S$ be a totally geodesic submanifold of a simply connected
    symmetric space $M.$ Then the following statements are
    equivalent:
    \begin{enumerate}
        \item $S$ is extrinsically symmetric;
        \item $S$ is adapted;
        \item $S$ is reflective.
    \end{enumerate}
\end{prop}

\begin{lem}\label{extrinsic symmetric} Let $\xi$ be of extrinsically symmetric type.
Assume $\lambda(M,\xi)$ to be odd or $M$ to be simply connected.
Then the submanifolds $M^{\xi}_t,\; t\notin \Z\pi$ are extrinsically
symmetric.
\end{lem}
\begin{proof}
Let $q=\Exp_p(t\Ad_{\LieG}(k)\xi),\, k\in\LieK_0$ be a point in
$M^{\xi}_t.$ If $\lambda(M,\xi)$ is odd then
$P^+_k=\Exp_p(\Ad(k)\lieP_+)$ is a meridian and hence reflective. If
$M$ is simply connected $P^+_k$ is adapted and by Proposition
\ref{extr-sym-adapted-reflective} reflective. Moreover
$T_qP^+_k=N_qM^{\xi}_t$ and $N_qP^+_k=T_qM^{\xi}_t.$ Let $r$ be the
reflection of $M$ in $P^+_k$ and let $\tau_k$ be the reflection of
$\lieP$ along $N_{\Ad_{\LieG}(k)\xi}M^{\xi}\cong\Ad(k)\lieP_+$ then
$\mathrm{D}r(p)=\tau_k.$ Take a point
$q_1=\Exp_p(t\Ad_{\LieG}(k_1)\xi),\, k_1\in\LieK_0$ in $M^{\xi}_t.$
Since $M^{\xi}$ is extrinsically symmetric in the euclidean space
$\lieP,$ there is an element $k_2\in\LieK_0$ such
$\tau_k\left(\Ad_{\LieG}(k_1)\xi\right)=\Ad_{\LieG}(k_2)\xi.$ Hence
$r(q_1)= r\left(\Exp_p(t\Ad_{\LieG}(k_1)\xi)\right)
=\Exp_p\left(\tau_k\left(t\Ad_{\LieG}(k_1)\xi\right)\right)=
\Exp_p\left(t\Ad_{\LieG}(k_2)\xi\right)$ is an element of
$M^{\xi}_t.$ Thus $r$ leaves $M^{\xi}_t$ invariant and
$\mathrm{D}r(q)X=X$ if $X$ is normal to $M^{\xi}_t$ and
$\mathrm{D}r(q)X=-X$ if $X$ is tangent to $M^{\xi}_t,$ i.e.\ $r$ is
the extrinsic symmetry of $M^{\xi}_t$ at $q.$
\end{proof}

\begin{rem}
  Assume $M$ to be an adjoint space $M=\overline{M}$ and moreover $M$ to be an {\it inner symmetric space},
  i.e.\ $s_p$ is contained in $\LieK_0,$ then $X$ lies in  $M^{\xi}$ if and only if $-X$ lies in $M^{\xi}.$ Although all
  submanifold $M^{\xi}_t,\; t\in (0,\pi)$ are symmetric spaces
  associated with the same orthogonal symmetric Lie algebra, the centrosome
  (polar) $M^{\xi}_{\pi\over 2}$ is not isomorphic to the other
  submanifolds. In fact, let $X\in M^{\xi},$ then $\gamma_X(t)=\gamma_{-X}(t)$ if and only
  if $t\in\Z{\pi\over 2}.$ Examples of this phenomenon are the spindles in real
  projective spaces, but also the spindle associated with
  $M^{\xi}=\LieSO(n)$ in the adjoint space of $\LieSO(2n)/
  (\LieSO(n)\times\LieSO(n))$ if $n$ is even.
\end{rem}

\begin{cor}
    Let $\xi$ be of extrinsically symmetric type. Assume $\lambda(M,\xi)$
    to be odd or $M$ to be simply connected. Then centrioles of chains
    of spindles are reflective.
\end{cor}

\begin{cor}\label{spindle=extrinsic symmetric}
    Let $\xi$ be of extrinsically symmetric type. Assume $\lambda(M,\xi)$
to be odd or $M$ to be simply connected. Then the submanifolds
$M^{\xi}_{\left(2n+1\right){\pi\over 2}-t}$
    and $M^{\xi}_{\left(2n+1\right){\pi\over 2}+t}, \; n\in\Z$ are extrinsically
    isometric.
\end{cor}
\begin{proof}
The reflection in $M^{\xi}_{\left(2n+1\right){\pi\over 2}}$ maps
$M^{\xi}_{\left(2n+1\right){\pi\over 2}-t}$ onto
$M^{\xi}_{\left(2n+1\right){\pi\over 2}+t}.$
\end{proof}

\section{Examples}\label{section-examples}

In this section we calculate the spindle number of all pairs $({M},
\xi),$ where ${M}$ is an indecomposable connected, simply connected
classical symmetric spaces of compact type and $\xi$ is of
extrinsically symmetric type. We use Cartan's notation for the
different types of the symmetric spaces $M.$ The spaces ${M}$ are
denoted as quotient spaces, but the given groups might not act
effectively. The exact notations as quotients of the identity
components of the isometry groups can be found in \cite{Wa-Zi-93}.
The classification of symmetric spaces and their symmetric
$s$-orbits on the Lie algebra level is done in \cite{Ko-Na-64} and
the elements $\xi$ of extrinsically symmetric type are given
explicitly in the classical cases. The list of the corresponding
pairs $\left({M},M^{\xi}\right)$ is taken from the appendix of
\cite{Be-Co-Ol-03}. To shorten our notations we denote by $J_n$ the
matrix $\left(\begin{array}{cc}
        0 & -I_n \\
        I_n & 0 \\
\end{array}\right).$

\begin{rem}
The Veronese type imbeddings of projective spaces over the number
fields $\R,\; \C$ or $\H$  as symmetric $s$-orbits and their chains
of spindles show that any natural number occurs as spindle number
(see type A I, A II and $\lieA$).
\end{rem}


\subsection{The classical non-hermitian case}

\subsubsection*{Type A I} Let
    ${M}=\LieSU(p+q)/\LieSO(p+q),\; 1\leq p\leq q$ and let $\xi=i\cdot
    \left(\begin{array}{cc}
                        a\cdot I_p & 0 \\
                        0 & b\cdot I_q \\
    \end{array}\right)$
    with $a=-{q\over {p+q}}$ and $b={p\over {p+q}}.$
    Then $\exp(t\xi)=\left(\begin{array}{cc}
                        \exp(ita)\cdot I_p & 0 \\
                        0 & \exp(itb)\cdot I_q \\
    \end{array}\right)$ lies in $\LieSO(p+q)$ if and only if
    $ta$ and $tb´$ are elements of $\pi\Z.$ Hence the spindle number is
    the smallest positive integer $n$
    such that $n{p\over{p+q}}$ and $n{q\over{p+q}}$ are
    integers.\footnotemark[1]

\subsubsection*{Type A II}
    Let ${M}=\LieSU(2(p+q))/\LieSp(p+q),\; 1\leq p\leq q$
    and let
    $\xi=\left(\begin{array}{cc}
                                E & 0\\ 0& E
     \end{array}\right)$
    with
    $E =i\left(\begin{array}{cc}
                        a\cdot I_p & 0 \\
                        0 & b\cdot I_q \\
     \end{array}\right),\;
    a=-{q\over {p+q}}$ and $b={p\over {p+q}}.$
    Then $\exp(t\xi)=\left(\begin{array}{cc}
                        M& 0 \\
                        0 & M \\
    \end{array}\right)$ with
    $M=\left(\begin{array}{cc}
                        \exp(ita)\cdot I_p & 0 \\
                        0 & \exp(itb)\cdot I_q \\
    \end{array}\right)$ lies in $\LieSp(p+q)$ if and only if
    $ta$ and $tb´$ are elements of $\pi\Z.$ Hence the spindle number is
    the smallest positive integer $n$
    such that $n{p\over{p+q}}$ and $n{q\over{p+q}}$ are
    integers.\footnotemark[1]

\subsubsection*{Type A III}
    Let ${M}=\LieSU(2n)/\mathrm{S}(\LieU(n)\times \LieU(n))$ and
    let $\xi={i\over 2}\cdot\left(
                \begin{array}{cc}
                    0 & I_n \\
                    I_n & 0 \\
                \end{array}
                \right).$
    Then $\exp(t\xi)=\cos\left({t\over 2}\right)\cdot
    \Id_{2n}+i\sin\left({t\over 2}\right)\cdot
    \left(\begin{array}{cc}
                        0& \Id_n \\
                        \Id_n & 0 \\
    \end{array}\right)$ lies in $\mathrm{S}(\LieU(n)\times \LieU(n))$ if
    and only if $t\in 2\pi\Z.$  Hence the spindle number is $2.$

\subsubsection*{Type BD I} We have to
distinguish two cases:
\begin{enumerate}
\item Let ${M}=\LieSO(p+q)/\LieSO(p)\times \LieSO(q),\; 1\leq p\leq q $
    and let
    $\xi= \left(\begin{array}{cc}
                    0 & E \\
                    -E^T & 0
     \end{array}\right)$ where $E$ is the rank one matrix
    $E=\left(\begin{array}{cc}
                    1 & 0 \\
                    0 & 0
     \end{array}\right)\in\mathrm{M}(p\times q).$
    An easy calculation shows that  $\exp(t\xi)$
     lies in $\LieSO(p)\times \LieSO(q) $ if and only if $\sin(t)=0$ and $\cos(t)=1.$  Thus
    $\lambda({M}, M^{\xi})=2.$
\item Let ${M}=\LieSO(2n)/\LieSO(n)\times \LieSO(n)$ and
    let
    $\xi={1\over 2}\cdot J_n.$
    Then $\exp(t\xi)= \cos\left({t\over 2}\right)\cdot
    \Id_{2n}+\sin\left({t\over 2}\right)\cdot J_n$
     lies in $\LieSO(n)\times \LieSO(n)$ if and only if
     $\sin\left({t\over 2}\right)=0$ and $\cos\left({t\over 2}\right)\cdot
    \Id_{2n}\in \LieSO(n)\times \LieSO(n).$  Thus
    $\lambda({M}, \xi)=2$ if $n$ is even and $\lambda(M,
    \xi)=4$ if $n$ is odd.
\end{enumerate}

\subsubsection*{Type D III}
    Let ${M}=\LieSO(4n)/\LieU(2n)$ and
    $\xi={1\over 2}\cdot\left(
            \begin{array}{cc}
                    J_n & 0 \\
                    0 & -J_n \\
            \end{array}
                \right).$
    Then $\exp(t\xi)=\cos\left({t\over 2}\right)\cdot
    \Id_{4n}+2\sin\left({t\over 2}\right)\cdot \xi$
    lies in $\LieU(2n)\cong\LieSO(4n)\cap\LieSp(2n) $ if and only if $t\in
    2\pi\Z.$  Hence the spindle number is $2.$

\subsubsection*{Type C I}
    Let ${M}=\LieSp(n)/\LieU(n)$ and let
    $\xi={i\over 2}\cdot \left(
            \begin{array}{cc}
                 I_n & 0 \\
                 0 & -I_n\\
            \end{array}
            \right).$
    Then $\exp(t\xi)=\cos\left({t\over 2}\right)\cdot
    \Id_{2n}+2\sin\left({t\over 2}\right)\cdot \xi$
    lies in $\LieU(n)\cong\LieSp(n)\cap\LieSO(2n) $ if and only if $t\in
    2\pi\Z.$  Hence the spindle number is $2.$

\subsubsection*{Type C II}
    Let ${M}=\LieSp(2n)/\LieSp(n)\times \LieSp(n)$ and
    $\xi={i\over 2}\cdot\left(
                \begin{array}{cccc}
                    0 & 0 & 0 & I_n \\
                     0 & 0 & I_n & 0 \\
                     0 & I_n & 0 & 0 \\
                    I_n & 0 & 0 & 0
                \end{array}
                \right).$
    Then $\exp(t\xi)=\cos\left({t\over 2}\right)\cdot
    \Id_{4n}+2\sin\left({t\over 2}\right)\cdot \xi$
    lies in $\LieSp(n)\times \LieSp(n)$ if and only if $t\in
    2\pi\Z.$  Hence the spindle number is $2.$


\subsection{The classical hermitian case}

Hermitian symmetric $s$-orbits arise as iso\-tropy orbits of simple
Lie groups:  A compact, connected, simply connected, simple Lie
group $\LieG$ endowed with a bi-invariant metric is a symmetric
space. The geodesic symmetry  $s_e$ at the identity $e$ is just the
inversion. The product $\LieG\times\LieG$ acts on $\LieG$ by
isometries: $(\LieG\times\LieG)\times \LieG\rightarrow \LieG,\,
(g,h)p=gph^{-1}.$ The involution on $\LieG\times \LieG$ given by the
conjugation with $s_e$ interchanges both factors. Thus its fixed
point set is the diagonal $\Delta^+\LieG$ in $\LieG\times\LieG.$ The
Cartan decomposition of the symmetric Lie algebra associated with
the symmetric space $\LieG\cong \LieG\times\LieG/\Delta^+\LieG$ is
$\lieG\times\lieG=\Delta^+\lieG\oplus \Delta^-\lieG,$ where
$\Delta^+\lieG$ denotes the diagonal and $\Delta^-\lieG$ the
antidiagonal in $\lieG\times\lieG.$ Let $\xi$ be an element of
$\lieG$ of extrinsically symmetric type then $(\xi, -\xi)$ is also
of extrinsically symmetric type in $\lieG\times\lieG.$  Moreover the
vector $(\xi,-\xi)$ defines a hermitian structure on its symmetric
$s$-orbit $M^{\xi}.$ To determine the spindle number of $(\LieG,
(\xi,-\xi))$ we have to look for the smallest positive integer $n$
such that $\exp(n\pi\cdot\xi)= \exp(-n\pi\cdot\xi).$

\subsubsection*{Type $\lieA$} Let $\LieG=\LieSU(p+q),\; 1\leq p\leq q$
and let $\xi=i\cdot \left(\begin{array}{cc}
                        a\cdot I_p & 0 \\
                        0 & b\cdot I_q \\
                \end{array}\right)$ with  $a=-{q\over {p+q}}$ and $b={p\over {p+q}}.$
Thus $\exp(t\xi)=\exp(-t\xi)$ if and only if $\exp(ita)=\exp(-ita)$
and $\exp(itb)=\exp(-itb).$ Hence the spindle number is the smallest
positive integer $n$ such that $n{p\over{p+q}}$ and $n{q\over{p+q}}$
are integers.\footnote[1]{Let $d$ be the greatest common divisor of
$p$ and
    $q,$ i.e.\ $p=dp_0$ and $q=dq_0,$ then $\lambda({M}, \xi)=p_0+q_0.$}

\subsubsection*{Type $\lieB$ and $\lieD$} Let $\LieG=\LieSpin(n).$ It is well known that $\LieSpin(n)$
is a two-fold covering of $\LieSO(n).$ The corresponding
identification of $\lieSpin(n)$ with $\lieO(n)$ described in
\cite{Frie-00} shows that the spindle number of $\LieSpin(n)$ with
respect to $\xi$ is two times the spindle number of $\LieSO(n)$ with
respect to $\xi.$ In the following we calculate the spindle number
of $\LieSO(n)$ with respect to $\xi= \left(\begin{array}{cc}
      0 & E \\
      -E^T & 0
\end{array}\right)$ where $E$ is the rank one matrix
$E=\left(\begin{array}{cc}
      1 & 0 \\
      0 & 0
\end{array}\right)\in\mathrm{M}(p\times q)$ with $p+q=n.$ Thus a straightforward calculation shows that
$\exp(t\xi)=\exp(-t\xi)$ if and only if $\sin(t)=0.$ Hence the
spindle number of $\LieSO(n)$ with respect to $\xi$ is $1$ and the
spindle number of $\LieSpin(n)$ with respect to $\xi$ is $2.$

\subsubsection*{Type $\lieC$} Let $\LieG=\LieSp(n)$
and $\xi={i\over 2}\cdot \left(\begin{array}{cc}
                                    I_n & 0 \\
                                    0 & -I_n\\
                         \end{array}\right).$
Thus $\exp(t\xi)=\cos\left({t\over 2}\right)\cdot \Id_{2n}
+2\sin\left({t\over 2}\right)\cdot \xi=\exp(-t\xi)$ if and only if
$\sin\left({t\over 2}\right)=0.$ Hence the spindle number is $2.$

\subsubsection*{Type $\lieD$} Let $\LieG=\LieSpin(2n).$
We determine the spindle of $\LieSO(2n)$ with respect to
$\xi={1\over 2}\cdot J_n.$ Thus $\exp(t\xi)=\cos\left({t\over
2}\right)\cdot \Id_{2n} +\sin\left({t\over 2}\right)\cdot
J_{n}=\exp(-t\xi)$ if and only if $\sin\left({t\over 2}\right)=0.$
Hence the spindle number of $\LieSO(2n)$ with respect to $\xi$ is
$2$ and the spindle number of $\LieSpin(2n)$ with respect to $\xi$
is $4.$

{\footnotesize
\begin{tabular}{|c|c|c|}
    \hline\rule[-2mm]{0mm}{7mm}
        ${\bf M}$
        &$\mathbf{M^{\xi}}$
        &{\bf spindle number}
        \\
    \hline
    \hline
        $\LieSU(p+q)/\LieSO(p+q)$
        &$\LieSO(p+q)/\mathrm{S}(\LieO(p)\times\LieO(q))$
        &$\begin{array}{l}
            \mathrm{smallest}\,n\in\N^*\\
            \mathrm{s.t.}\; n{p\over{p+q}}\in\N^*\\
            \mathrm{and}\;n{q\over{p+q}}\in\N^*\footnotemark[1]\\
        \end{array}$
        \\
    \hline
        $\LieSU(2(p+q))/\LieSp(p+q)$
        &$\LieSp(p+q)/\LieSp(p)\times\LieSp(q)$
        &$\begin{array}{l}
            \mathrm{smallest}\,n\in\N^*\\
            \mathrm{s.t.}\; n{p\over{p+q}}\in\N^*\\
            \mathrm{and}\;n{q\over{p+q}}\in\N^*\footnotemark[1]\\
        \end{array}$
        \\
    \hline
        $\LieSU(2n)/\mathrm{S}(\LieU(n)\times \LieU(n)) $
        &$\LieU(n)$
        &2
        \\
   \hline
        $\LieSO(p+q)/\LieSO(p)\times \LieSO(q) $
        &$(S^{p-1}\times S^{q-1})/\Z_2$
        & 2
        \\
    \hline
        $\LieSO(2n)/\LieSO(n)\times \LieSO(n)$
        &$\LieSO(n)$
        &$\begin{array}{l}
            2,\; \mathrm{if}\; n \;\mathrm{ is\; even}\\
            4,\; \mathrm{if}\; n \;\mathrm{ is\; odd}\\
        \end{array}$
        \\
    \hline
        $\LieSO(4n)/\LieU(2n) $
        &$\LieU(2n)/\LieSp(n)$
        &2
        \\
    \hline
        $\LieSp(n)/\LieU(n) $
        &$\LieU(n)/\LieSO(n)$
        &2
        \\
    \hline
        $\LieSp(2n)/\LieSp(n)\times \LieSp(n) $
        &$\LieSp(n)$
        &2
        \\
    \hline
    \hline
        $\LieSU(p+q)$
        &$\LieSU(p+q)/\mathrm{S}(\LieU(p)\times\LieU(q))$
        &$\begin{array}{l}
            \mathrm{smallest}\,n\in\N^*\\
            \mathrm{s.t.}\; n{p\over{p+q}}\in\N^*\\
            \mathrm{and}\;n{q\over{p+q}}\in\N^*\footnotemark[1]\\
        \end{array}$\\
    \hline
        $\LieSpin(n)$
        &$\LieSO(n)/(\LieSO(2)\times \LieSO(n-2))$
        &2\\
    \hline
        $\LieSp(n) $
        &$\LieSp(n)/\LieU(n)$
        &2\\
    \hline
        $\LieSpin(2n)$
        &$\LieSO(2n)/\LieU(n)$
        &4\\
    \hline
\end{tabular}}

\end{document}